\documentclass[11pt,a4paper]{article}
\usepackage[T1]{fontenc}
\usepackage[utf8]{inputenc}
\usepackage{lmodern}
\usepackage{microtype}
\usepackage[margin=1.12in]{geometry}
\usepackage{amsmath,amssymb,amsthm,mathtools}
\usepackage{mathrsfs}
\usepackage{enumitem}
\usepackage{xcolor}
\usepackage{booktabs}
\usepackage{array}
\usepackage[most]{tcolorbox}
\usepackage{aliascnt}
\usepackage{hyperref}
\usepackage[capitalize,nameinlink,noabbrev]{cleveref}
\hypersetup{colorlinks=true,linkcolor=blue!40!black,citecolor=blue!40!black,urlcolor=blue!40!black,pdftitle={Annular Obstructions and Skew-Shape Block Identities},pdfauthor={Adam Tilepe}}
\setlist{nosep,topsep=4pt}
\numberwithin{equation}{section}
\newtheorem{theorem}{Theorem}[section]
\newaliascnt{lemma}{theorem}
\newtheorem{lemma}[lemma]{Lemma}
\aliascntresetthe{lemma}
\newaliascnt{proposition}{theorem}
\newtheorem{proposition}[proposition]{Proposition}
\aliascntresetthe{proposition}
\newaliascnt{corollary}{theorem}
\newtheorem{corollary}[corollary]{Corollary}
\aliascntresetthe{corollary}
\newaliascnt{definition}{theorem}
\newtheorem{definition}[definition]{Definition}
\aliascntresetthe{definition}
\newaliascnt{remark}{theorem}
\newtheorem{remark}[remark]{Remark}
\aliascntresetthe{remark}
\newtcolorbox{majorresult}{
  enhanced,breakable,
  colback=blue!1!white,
  colframe=blue!32!black,
  boxrule=0.5pt,
  borderline west={2.0pt}{0pt}{blue!60!black},
  arc=1.2pt,
  left=9pt,right=9pt,top=7pt,bottom=7pt,
  before skip=11pt,after skip=11pt,
  parbox=false
}

\newcommand{\ZZ}{\mathbb Z}
\newcommand{\QQ}{\mathbb Q}

\newcommand{\JJ}{\mathcal J}
\newcommand{\HH}{\mathcal H}
\newcommand{\RR}{\mathcal R}
\newcommand{\Om}{\Omega}
\newcommand{\bl}{\widetilde{\mathrm{bl}}}

\newcommand{\Conn}{\operatorname{Conn}}

\title{\bfseries Annular Obstructions and Skew-Shape Block Identities\\[3pt]
\large Weighted Order Ideals and a Lean-Verified Refinement}
\author{Adam Til\`epe}
\date{11 November 2026}
\begin{document}
\maketitle
\begin{abstract}
We give three principal results on weighted order ideals, cluster
coefficients, and block-statistic enumerators. First, in the rank-two
annular exchange model, the frozen row $(2,-3)$ yields the specialized
polynomial $1+q^2$, which is not unimodal. Its interpretation as a
single-lamination counterexample to Theorem~14 of Kang, Lee and Lim
(\emph{European Journal of Combinatorics}, 2026) depends on identifying
this row with the shear coordinates of a bridging curve in the same
triangulation. We discuss this geometric correspondence using published
annular examples; it is not established by the algebraic calculation or
formalized in Lean. Second, the same polynomial fails log-concavity,
and the frozen row $(1,-2)$ gives the full-support polynomial
$1+2q+7q^2+10q^3+9q^4+4q^5+q^6$, likewise violating log-concavity.
These calculations provide obstructions relevant to their Conjecture~1,
subject to the same geometric identification. Third, for arbitrary
finite skew Young diagrams, we prove the block-statistic identity of
Kahane's Conjecture~4.7 (arXiv:2607.11225v2), using an explicit
correction of the root convention. The proof uses commuting
tableau-transfer matrices, a connected-block inverse kernel, and
finite polynomial differential equations, with companion Lean~4/Mathlib
formalizations. As a consequence, we obtain a bivariate block
refinement that is coefficientwise nonnegative in $t-1$ after a common
monomial shift in $q$.
\end{abstract}

\noindent\textbf{Keywords.} Order polynomials; skew Young diagrams; order ideals;
block statistics; cluster algebras from surfaces; log-concavity;
dual stable Grothendieck polynomials; formal proof.\\
\textbf{MSC 2020.} 05A15, 06A07, 05E05, 13F60.

\section{Introduction}\label{sec:intro}
For a finite poset $P$, write $\JJ(P)$ for its set of lower order
ideals and $\Om(P;t)$ for its order polynomial, characterized by the
number of weakly order-preserving maps $P\to\{1,\ldots,m\}$ at each
positive integer $m$.  The ordinary rank polynomial and the order
polynomial both arise from $\JJ(P)$, but record different information.
Replacing cardinality by a general integer weight on the elements of
$P$ makes that difference particularly sharp.

Our first point of contact is the order-ideal expansion of cluster
variables on marked surfaces.  Its principal-coefficient specialization
produces an ordinary rank polynomial, whereas a geometric coefficient
corresponding to a lamination may give nonuniform weights and internal
zeros.  Kang, Lee and Lim \cite{KangLeeLim2026} assert unimodality for
single-lamination specialized cluster expansions and conjecture
log-concavity.  We explain why a bridging curve in an annulus permits an
initial coefficient of magnitude two and compute two explicit
obstructions to these assertions.  The geometric identification uses
the negative-shear coordinates in Felikson and Tumarkin
\cite[Fig.~4.2]{FeliksonTumarkin2024}; this sign conversion is
spelled out in \cref{sec:annulus}.

Our second point of contact is Kahane's construction of an intrinsic
block statistic on labeled fence posets \cite{Kahane2026}.  Its
extension to skew-shape cell posets is formulated as
Conjecture~4.7 in the second version of Kahane's paper
(arXiv:2607.11225v2).
The $2\times n$ rectangle case was established earlier by White
\cite{White2026}. We give a transfer-matrix treatment for arbitrary
skew shapes using the root-exclusion convention of \cref{def:blocks},
corresponding to Kahane's Definition~4.5 in v2; see
\cref{sec:root-convention,sec:scope}.

The two constructions are joined by one bivariate enumerator:
\[
  \HH_{P,w}(t;q)
  =\sum_{I\in\JJ(P)}q^{w(I)}\Om(P\setminus I;t-1).
\]
At $q=1$ it recovers the order polynomial; at $t=2$ it recovers
the weighted ideal polynomial.  The skew-shape block formula also yields, as a direct consequence,
a positive expansion of \eqref{eq:intro-H}, but \emph{does not} restore
unimodality of arbitrary weighted specializations.  The one-cell case
already displays the distinction: $(t-1)+q^2$.

The paper follows the common poset through weighted ideals
(\cref{sec:ideals}), annular examples (\cref{sec:annulus}),
canonical blocks (\cref{sec:blocks}), transfer matrices
(\cref{sec:transfer}), and the resulting block identity
(\cref{sec:main}).  \Cref{sec:scope} records the mathematical scope,
and \cref{sec:formalization} documents the companion Lean development.

\section{Weighted ideals and order polynomials}\label{sec:ideals}
Let $P$ be a finite poset and $w:P\to\ZZ$.  For $I\in\JJ(P)$ put
$w(I)=\sum_{p\in I}w(p)$ and define
\begin{equation}\label{eq:weighted-rank}
  \RR_{P,w}(q)=\sum_{I\in\JJ(P)}q^{w(I)}\in\ZZ[q,q^{-1}].
\end{equation}
The usual rank polynomial is $\RR_{P,\mathbf1}(q)$.
We now set
\begin{equation}\label{eq:intro-H}
  \HH_{P,w}(t;q)
  =\sum_{I\in\JJ(P)}q^{w(I)}\Om(P\setminus I;t-1).
\end{equation}
We use $\Om(\varnothing;t)=1$ and $\Om(Q;0)=0$ for nonempty $Q$.
Equation \eqref{eq:intro-H} defines an element of
$\QQ[t,q,q^{-1}]$; for an evaluation with a negative exponent of $q$
we assume $q\ne0$.

\begin{majorresult}
\begin{proposition}[The two specializations]\label{prop:slices}
For every finite poset $P$ and integer weighting $w$,
\begin{equation}\label{eq:slices}
  \HH_{P,w}(t;1)=\Om(P;t),\qquad
  \HH_{P,w}(2;q)=\RR_{P,w}(q).
\end{equation}
In particular $\RR_{P,\mathbf1}(1)=\Om(P;2)=|\JJ(P)|$.
\end{proposition}
\end{majorresult}
\begin{proof}
For each positive integer $m$, a weakly order-preserving map
$f:P\to[m]$ determines the ideal $I=f^{-1}(1)$; on $P\setminus I$
the function $f-1$ is an order-preserving map to $[m-1]$.  Conversely
these data reconstruct $f$.  Consequently
$\Om(P;m)=\sum_{I\in\JJ(P)}\Om(P\setminus I;m-1)$.
The identity holds for infinitely many $m$, so it holds as an identity
of polynomials.  This gives the first specialization.  The second
follows from $\Om(Q;1)=1$ for every finite $Q$.
\end{proof}

Write $Z$ for the zeta matrix of the lattice $\JJ(P)$:
$Z[I,J]=\mathbf1_{I\subseteq J}$.
\begin{lemma}[Weak ideal chains]\label{lem:zeta}
For every positive integer $m$,
\begin{equation}\label{eq:zeta}
  \Om(P;m)=(Z^m)[\varnothing,P].
\end{equation}
\end{lemma}
\begin{proof}
Matrix multiplication counts chains
$\varnothing=I_0\subseteq I_1\subseteq\cdots\subseteq I_m=P$.
Such chains are in bijection with isotone maps $P\to[m]$ by recording
the stage at which each element first enters the chain.
\end{proof}

\subsection{Fences and skew ribbons}
A \emph{fence} is a poset whose Hasse diagram is an oriented path.
A \emph{skew ribbon} is a connected skew Young diagram with no
$2\times2$ square.  Its cell order is the coordinatewise order,
with one fixed choice of drawing convention.
\begin{lemma}[Fence--ribbon correspondence]\label{lem:fence}
Every finite fence poset admits a skew-ribbon realization, and the
cell poset of any skew ribbon has an oriented-path Hasse diagram.
\end{lemma}
\begin{proof}
Enumerate the elements in their path order as $p_1,\ldots,p_n$.
Place the first cell at $(0,0)$ in a grid whose second coordinate
increases eastward and whose first coordinate increases southward.
For each cover $p_i<p_{i+1}$, append the next cell one step east;
for each cover $p_i>p_{i+1}$, append it one step north.  The resulting
sequence of cells is a northeast-moving edge-connected path.  Its
occupied cells form row intervals with boundaries weakly decreasing
from top to bottom, so after translation it is a skew Young diagram.
Because only north and east steps occur, it contains no $2\times2$
square.  Consecutive cells are comparable, with the required cover
orientation; any nonconsecutive cells are either incomparable or have
an intermediate cell on the path, so cannot form an extra cover.
Thus its Hasse diagram is precisely the prescribed oriented path.
Conversely, a connected skew ribbon can be traced from its southwest
to its northeast end using north and east steps.  Every edge of its
Hasse diagram occurs between consecutive cells of this trace, and no
other cover is possible; the trace therefore recovers a fence.
\end{proof}

\begin{lemma}[No gaps for small weights]\label{lem:no-gaps}
If $|w(p)|\le1$ for every $p\in P$, every integer between
$\min_{I\in\JJ(P)}w(I)$ and $\max_{I\in\JJ(P)}w(I)$ is the weight of
an ideal.  Thus the normalized weighted rank polynomial has no
internal zeros.
\end{lemma}
\begin{proof}
The undirected cover graph of $\JJ(P)$ is connected.  Choose a path
from an ideal of minimum weight to one of maximum weight.  Two ideals
joined by a cover differ by a single element, hence the weight changes
by at most one.  Such a walk must visit every intermediate integer.
\end{proof}

\begin{remark}\label{rem:no-gaps}
The conclusion is weaker than unimodality.  For a six-element chain,
the successive weights $(-1,1,1,1,1,-1)$ give
$q\RR_{P,w}(q)=1+2q+q^2+2q^3+q^4$, whose coefficients are not
unimodal despite having full support.  This is an abstract weighted
poset example; no geometric lamination is asserted here.
\end{remark}

\subsection{Surface coefficients as weighted ideals}
For a plain arc $\gamma$ admitting a principal-coefficient expansion
of the form
\begin{equation}\label{eq:F}
 F_\gamma(y_1,\ldots,y_s)
  =\sum_{I\in\JJ(P_\gamma)}\prod_{p\in I}y_{\tau(p)},
\end{equation}
assign an integer shear vector $(s_1,\ldots,s_s)$ and set
$w(p)=s_{\tau(p)}$.  The substitution $y_i=q^{s_i}$ yields
$\RR_{P_\gamma,w}(q)$.  The separation formula of Fomin and
Zelevinsky \cite[Thm.~3.7]{FominZelevinsky2007} explains the
normalizing tropical monomial in the corresponding cluster expansion.
For the normalization exponent $m_\gamma=\min_{I}w(I)$, one obtains
\begin{equation}\label{eq:surface-slice}
 c_{\gamma,L}(q)=q^{-m_\gamma}\RR_{P_\gamma,w}(q)
                =q^{-m_\gamma}\HH_{P_\gamma,w}(2;q)
\end{equation}
whenever the particular geometric arc and lamination have been
identified with the expansion \eqref{eq:F}.  The existence of this
identification is not a consequence of the algebraic identity alone;
see also the surface formulas of Musiker, Schiffler and Williams
\cite{MusikerSchifflerWilliams2011}.

\section{Annular obstructions to unimodality and log-concavity}
\label{sec:annulus}
\begin{theorem}[Kang--Lee--Lim, Theorem~14; original assertion]
\label{thm:published}
For a single lamination on an unpunctured marked surface, the
associated specialized cluster $c$-polynomial is unimodal.
\end{theorem}
Kang, Lee and Lim also propose the log-concavity of the same
coefficient sequence as their \emph{Conjecture~1}
\cite[\S5.1]{KangLeeLim2026}.  We keep the numbering of their original
statements to make the comparison unambiguous.

\subsection{The induction base and the annular geometry}
In the arXiv v3 proof of \cref{thm:published}, the initial polynomial
$c_{x_1}$ is asserted to be $1$ or $1+q$.  This is not implied by the
unrestricted hypothesis of a single lamination: shear coordinates of
an individual curve need not lie in $\{-1,0,1\}$.

Consider the annulus $S_{1,1}$ with one marked point on each boundary
component.  It has an ideal triangulation with principal exchange matrix
\begin{equation}\label{eq:kronecker}
 B=\begin{pmatrix}0&2\\-2&0\end{pmatrix}.
\end{equation}
Felikson and Tumarkin \cite[Fig.~4.2 and Lemma~4.1]{FeliksonTumarkin2024}
draw the annulus, its universal cover, and individual bridging curves
with coefficient vectors
$(-2,3),(-1,2),(0,1),\ldots$.  Their vectors are explicitly
\emph{negative} shear coordinates.  After converting to the
shear convention used in the frozen rows of
\cite[Defs.~11--12]{KangLeeLim2026}, the corresponding rows are
$(2,-3),(1,-2),(0,-1),\ldots$, up to the matching choice of arc
orientation.  More precisely, Felikson and Tumarkin state in the proof of
\cite[Lemma~4.1]{FeliksonTumarkin2024} that every bridging curve on
$S_{1,1}$ is obtained from a fixed one by an \emph{integer} power
of the Dehn twist along the core curve. Both positive and negative
powers are allowed, not merely the finitely many twists drawn in
Fig.~4.2. Following the branch starting with the displayed negative
shear vectors $(0,1),(-1,2),(-2,3)$ gives
$(-r,r+1)$ for $r\ge0$, and hence the candidate rows
$(r,-(r+1))$ after reversing signs. This use of the twist orbit
concerns individual bridging curves, not admissibility of
mutation-finite coefficient vectors.

We checked these conventions with the authors in correspondence on
10 November 2026. A.~Felikson confirmed that one may reverse the signs
of the displayed shear vectors while preserving their coordinate order,
and that each depicted curve is individually a lamination. W.~Kang
confirmed that the surface setting in \cite{KangLeeLim2026} has no
topological restriction beyond being unpunctured, while emphasizing that
the shear coordinates of a \emph{single} lamination satisfy an alternating
condition. He also confirmed that Definition~20 obtains the
$c$-polynomial by specializing mutable cluster variables to $1$ and
the unique frozen variable to $q$, without additional normalization.
These clarifications concern the stated conventions, not an independent
verification of the counterexamples. The matrix-mutation calculations
remain algebraic and are checked in the companion Lean development.

For the algebraic calculation, the proposed extended exchange matrix is
\begin{equation}\label{eq:extended}
 \widetilde B_r=
 \begin{pmatrix}0&2\\-2&0\\r&-(r+1)\end{pmatrix}.
\end{equation}
The two off-diagonal entries $\pm2$ of its principal part record the
two ideal triangles of the annular triangulation. The third row is the
proposed shear-coordinate row of one bridging curve after the sign
conversion described above. Kang--Lee--Lim define a single lamination in
\cite[Def.~10]{KangLeeLim2026}, specify shear coordinates by the
local $S$- and $Z$-shaped crossing contributions in
\cite[Def.~11]{KangLeeLim2026}, and place these coordinates in the
frozen rows of the extended matrix in
\cite[Def.~12]{KangLeeLim2026}. These definitions do not formulate
a separate two-entry ``alternating-sign condition'' that could be
verified merely by checking opposite signs. For $r\ge1$ the proposed
row has opposite signs, while for $r=0$ one entry is zero; neither
observation replaces the local shear-coordinate calculation.
Kang emphasized an alternation constraint in correspondence, but
its precise applicability to this indexed row is not inferred from
the sign pattern alone. To interpret the following algebraic calculation
as a surface counterexample, the row, its coordinate order, and the
principal matrix must all refer to the \emph{same} triangulation.
This geometric identification is an additional step, not established
by the formal matrix-mutation check.

\begin{majorresult}
\begin{theorem}[Annular non-unimodality]\label{thm:annulus}
In the rank-two exchange model \eqref{eq:kronecker}, take the
frozen coefficient row $(r,-r-1)$, with $r\ge2$.
The first mutation produces
\begin{equation}\label{eq:1plusqr}
 c_r(q)=1+q^r.
\end{equation}
In particular the coefficient sequence is not unimodal.  For $r=2$
it also fails log-concavity. If the row is identified with the
shear coordinates of a single bridging lamination in the same
triangulation, these yield counterexamples to \cref{thm:published}
and Kang--Lee--Lim Conjecture~1 in that geometric realization.
\end{theorem}
\end{majorresult}
\begin{proof}
For the extended exchange matrix \eqref{eq:extended},
the exchange relation at the first mutable arc is
$x_1'=(x_2^2+q^r)/x_1$.  Setting $x_1=x_2=1$ gives
\eqref{eq:1plusqr}.  The coefficient sequence has positive
coefficients in degrees $0$ and $r$ and zero coefficients in between;
it is not unimodal when $r\ge2$.  At $r=2$ the central inequality is
$a_1^2=0<a_0a_2=1$.
\end{proof}

\begin{remark}
For $r>2$, non-unimodality does not by itself force a failure of
the three-term log-concavity inequality: the direct counterexample
to that inequality is $r=2$.  In contrast, the next example has no
internal zero coefficients at all.
\end{remark}

\begin{majorresult}
\begin{proposition}[Full-support failure of log-concavity]
\label{prop:positive}
For the single-curve frozen row $(1,-2)$, the rank-two mutation
sequence $(1,2,1,2)$ gives
\begin{equation}\label{eq:full-support}
 c(q)=1+2q+7q^2+10q^3+9q^4+4q^5+q^6.
\end{equation}
All coefficients are strictly positive, but $a_1^2=4<7=a_0a_2$.
\end{proposition}
\end{majorresult}
\begin{proof}
Write $\widetilde B^{(j)}$ for the extended exchange matrix
after the first $j$ mutations, starting with
$\widetilde B^{(0)}=\left(\begin{smallmatrix}0&2\\-2&0\\1&-2\end{smallmatrix}\right)$.
Applying the ordinary extended-matrix mutation rule in sequence
$1,2,1,2$ gives the following complete audit of the frozen row;
the mutable $2\times2$ part alternates sign at each step:
\[
\renewcommand{\arraystretch}{1.12}
\begin{array}{c|c|c|c}
\text{step }j&\text{frozen row before step}&
  (b_{12},b_{21})\text{ before step}&
  \text{frozen row after step}\\ \hline
1&(1,-2)&(2,-2)&(-1,0)\\
2&(-1,0)&(-2,2)&(-1,0)\\
3&(-1,0)&(2,-2)&(1,0)\\
4&(1,0)&(-2,2)&(1,0)
\end{array}
\]
The column signs in these rows yield, \emph{before setting the
mutable variables to $1$}, precisely the four exchange equations
\[
\begin{aligned}
 x_1x'_1&=q+x_2^2,&
 x_2x'_2&=1+(x'_1)^2,\\
 x'_1x''_1&=1+q(x'_2)^2,&
 x'_2x''_2&=1+(x''_1)^2.
\end{aligned}
\]
By the Laurent phenomenon, every mutated cluster variable belongs
to $\ZZ[q^{\pm1}][x_1^{\pm1},x_2^{\pm1}]$, so evaluation at
$x_1=x_2=1$ is a well-defined ring homomorphism. The exact divisions
below are also checked directly in $\ZZ[q]$ by multiplication; the
Laurent phenomenon alone is not being used to assert polynomiality in
$q$. Set the initial $x_1=x_2=1$, only after forming these rational
expressions, and denote the four successive values by
$f_1,f_2,f_3,f_4$.  The exchange equations become
\begin{align*}
 f_1&=1+q,\\
 f_2&=1+f_1^2=2+2q+q^2,\\
 f_3&=\frac{1+q f_2^2}{f_1}
       =1+3q+5q^2+3q^3+q^4,\\
 f_4&=\frac{1+f_3^2}{f_2}
       =1+2q+7q^2+10q^3+9q^4+4q^5+q^6.
\end{align*}
At step~3 the active-column frozen entry is $-1$, which places
$q^1$ beside $(x'_2)^2$ in the numerator of the exchange relation;
at step~4 the active-column frozen entry is $0$, so no power of
$q$ appears in that exchange relation. Both divisions are exact as verified by
multiplying out their displayed polynomial identities.  Finally,
$[q]f_4=2$ and $[q^2]f_4=7$, so
$([q]f_4)^2=4<[q^0]f_4\,[q^2]f_4=7$.
\end{proof}

\begin{majorresult}
\begin{corollary}[Scope of the published assertions]\label{cor:refute}
Assume that the frozen rows in \cref{eq:extended} are identified
with the shear coordinates of individual bridging laminations, with
the stated order of arcs and the principal matrix belonging to the
same annular triangulation. Under this geometric identification,
\cref{thm:published} and the log-concavity assertion of
Kang--Lee--Lim Conjecture~1 both admit counterexamples.
At the algebraic level, the first mutated polynomial $1+q^2$
shows why the induction-base restriction used in the printed proof
of Theorem~14 does not follow from the coefficient-row calculation
alone.
\end{corollary}
\end{majorresult}
\begin{proof}
Apply \cref{thm:annulus,prop:positive} under the stated
geometric identification. The annulus is unpunctured, and each
proposed lamination consists of one bridging curve. The
initial-coefficient calculation was given in
\cref{eq:1plusqr}.
\end{proof}

\begin{majorresult}
\begin{proposition}[The one-cell bridge]\label{prop:one-cell}
Let $P=\{p\}$ and $w(p)=2$.  Then
\begin{equation}\label{eq:one-cell}
 \HH_{P,w}(t;q)=(t-1)+q^2.
\end{equation}
Consequently $\HH_{P,w}(2;q)=1+q^2$ while
$\HH_{P,w}(t;1)=t=\Om(P;t)$.
\end{proposition}
\end{majorresult}
\begin{proof}
The two ideals are $\varnothing$ and $\{p\}$.  Their terms in
\eqref{eq:intro-H} are $\Om(P;t-1)=t-1$ and
$q^2\Om(\varnothing;t-1)=q^2$.
\end{proof}

\begin{remark}
These examples concern nonuniform \emph{geometric} weights.  They
do not contradict the results of \cite{KangLeeLim2026} about the
unweighted rank polynomials of ordinary fence or loop-fence posets.
Nor do they rule out a corrected coefficient-shape theorem for a
smaller geometric class, such as certain polygonal laminations.
\end{remark}

\section{Canonical blocks in skew diagrams}\label{sec:blocks}
Let $D=\lambda/\mu$ be a finite skew Young diagram, with cell poset
$P_D$ ordered by the chosen coordinatewise convention.  Disconnected
shapes are allowed.  For any skew subdiagram $E$, write $h(E)$ for its number of occupied
rows and $\operatorname{col}(E)$ for its number of occupied columns.  These
symbols are distinct from the integer weighting $w:P_D\to\ZZ$ and
from the exchange matrix $B$ used in \cref{sec:annulus}.

\subsection{The root convention}\label{sec:root-convention}
Kahane's skew-shape block statistic is defined in
\cite[Def.~4.5]{Kahane2026}, using the global block-order condition
of his Definition~3.2.  The printed strict boundary inequalities
quantify over every cell, and hence include the root itself.  We use
the mathematically consistent root-excluded interpretation, which
Kahane confirmed in correspondence on 10 November 2026; the precise
comparison with version~2 is recorded in \cref{sec:scope}.

\begin{definition}[Root-excluded block statistic]\label{def:blocks}
For a bijective labeling $\sigma:D\xrightarrow\sim[n]$, a block is
a connected skew subdiagram $B$ with its southwest root $r(B)$.
Every nonroot cell at the left end of an occupied block row has
label strictly below $\sigma(r(B))$; every nonroot cell at the
bottom end of an occupied block column has label strictly above it.
The order of distinct blocks is compatible with the increasing
order of their root labels.  The canonical number of blocks is
denoted $\bl_D(\sigma)$.
\end{definition}

The equivalent computational definition assigns a root label to each
cell as one scans from the southwest boundary toward the northeast.
The lower-neighbor root $a$ and the left-neighbor root $b$ constrain
the new root, and the median operation implements the convention.
This makes the definition precise even in an automated formal proof.

\begin{lemma}[Boundary cells of a connected skew shape]\label{lem:boundary}
Let $B$ be a nonempty connected skew diagram with $h$ occupied rows
and $\operatorname{col}(B)$ occupied columns.  It has exactly $h-1$ nonroot leftmost row
cells and $\operatorname{col}(B)-1$ nonroot bottommost column cells; these two sets are
disjoint.
\end{lemma}
\begin{proof}
Number the occupied rows from top to bottom and let their column
intervals be $[\ell_i,r_i]$.  For a skew Young diagram both endpoint
sequences are weakly decreasing down the rows; for a connected diagram,
consecutive intervals overlap.  If row $i$ is not the bottom row, then
$\ell_{i+1}\le\ell_i\le r_{i+1}$ by monotonicity and overlap.
Consequently the leftmost cell $(i,\ell_i)$ has a cell immediately
below it, so it cannot be a bottom-column cell.  The leftmost cell
of the bottom occupied row has no cell below and is simultaneously
a left-row boundary and a bottom-column boundary; it is the unique
intersection of these boundary sets.  There are $h$ left-row
boundaries and $\operatorname{col}(B)$ bottom-column boundaries, one for
each occupied row and column.  Excluding their unique common root
gives $h-1$ and $\operatorname{col}(B)-1$ cells in disjoint sets.
\end{proof}

\begin{lemma}[Median construction and uniqueness]\label{lem:median}
Under \cref{def:blocks}, every bijective labeling of a skew diagram
has a unique valid block partition.  For a cell with label $x$,
write $a$ for the root label at its lower neighbor and $b$ for that
at its left neighbor, using $-\infty$ or $+\infty$ when the
corresponding neighbor is absent.  Its canonical root label is
\begin{equation}\label{eq:median}
  \rho(c)=\operatorname{median}(a,x,b).
\end{equation}
A new root is created precisely when $a<x<b$.
\end{lemma}
\begin{proof}
Orient row indices downward and column indices rightward.  Write
$\rho(c)$ for the provisional root label at a processed cell.  We
maintain two invariants throughout the southwest-first scan:
\begin{equation*}
 \rho(\text{below }c)\le \rho(c)\le\rho(\text{left of }c)
 \quad\text{whenever the indicated neighbor exists},
\end{equation*}
and every occupied predecessor has already received a root label.
These assertions are vacuous before the first step.

Here is the crucial compatibility of the two bounds.  When $c$ has
both a lower neighbor $d$ and a left neighbor $\ell$, the southwest
corner $s$ of the unit square with opposite corner $c$ also belongs
to $D$: row intervals and decreasing left endpoints force its
presence.  Its root has already been assigned, and the local
inequalities along $s\to d$ and $s\to\ell$ yield
\[
 a=\rho(d)\le\rho(s)\le\rho(\ell)=b.
\]
Thus $a\le b$ follows from the preserved local order, not merely
from the shape of $D$.  The median $\rho(c)=\max\{a,\min\{x,b\}\}$
then satisfies $a\le\rho(c)\le b$.  With only a lower (respectively,
left) neighbor the rule becomes $\max(a,x)$ (respectively,
$\min(x,b)$); with neither neighbor the cell becomes a new root.
In each case the two local inequalities persist.  This proves the
invariants by induction over the scan.

Every cell either starts a block at its own label or inherits the
root of a previously processed lower or left neighbor.  Following
inherited edges therefore connects it to its anchor.  Moreover the
root labels are weakly monotone along every occupied row and column.
A fixed root's cells consequently form row intervals and column
intervals; the same monotonicity across comparable corners gives
weakly decreasing row endpoints.  Each nonempty fiber is thus a
connected skew subdiagram, and its southwest anchor is its unique
root.  The three median possibilities give the prescribed nonroot
boundary-label inequalities.  More explicitly, when $x<a$, the cell joins the
lower-neighbor fiber at the left end of its row, so $x<a$ is
the required left-boundary inequality.  When $x>b$, it joins
the left-neighbor fiber at the bottom end of its column, so
$x>b$ is the required bottom-boundary inequality.  When
$a<x<b$, it creates a new singleton root.  If $a=b$, both neighboring fibers coincide,
and the new cell joins them.  The local inequalities extend
along row and column segments, and through the present southwest
corner for two comparable cells, to the global root order required
in Kahane's Definition~3.2.

Conversely, suppose a partition satisfies the root-excluded boundary
inequalities and the global order.  The latter implies the weak local
inequalities for adjacent cells.  At an unprocessed cell, its root
must be either its own label or the root of an earlier neighboring
block (connectedness precludes any other newly appearing fiber).
The boundary inequalities force the lower root for $x<a$, the left
root for $x>b$, and a new root for $a<x<b$; when $a=b$ both
predecessors belong to the same root fiber.  Missing neighbors are
handled by the corresponding one-sided max/min rule.  Hence every
valid partition agrees with the median scan at each step, proving
uniqueness.  For a disconnected diagram the same argument applies to
each connected component.  The separate formal lemmas behind these
inductive steps are identified in \cref{sec:formalization}.
\end{proof}

\begin{remark}\label{rem:kahane-convention}
\Cref{lem:median} proves the explicitly root-excluded statistic.
The corresponding conjecture is Conjecture~4.7 of
arXiv:2607.11225v2, proved here as \cref{thm:kahane} for the corrected
statistic.  Its Definition~4.5 is compared with our coordinate
convention in \cref{sec:scope}.
\end{remark}

\section{Commuting tableau-transfer matrices}\label{sec:transfer}
All matrices in this section are indexed by the finite set
$\JJ(P_D)$.  For $I\subseteq J$ the difference $J\setminus I$ is
a skew diagram.  Use $\operatorname{col}(B)$ and $h(B)$ for the numbers of occupied columns
and rows, respectively, with $\operatorname{col}(\varnothing)=0$.  Define
\begin{align}
 T_x[I,J]&=\begin{cases}x^{\operatorname{col}(J\setminus I)},&I\subseteq J,\\
                           0,&\text{otherwise},\end{cases}
                   \label{eq:T}\\
 A_z[I,J]&=\begin{cases}
 z^{h(J\setminus I)-1}(1-z)^{\operatorname{col}(J\setminus I)-1},
   & I\subsetneq J,\ \Conn(J\setminus I),\\
 0,&\text{otherwise}.
 \end{cases}\label{eq:A}
\end{align}
In particular $T_0=\mathrm{Id}$ and $T_1=Z$, while every $A_z$
is strictly upper triangular with respect to inclusion of ideals.

\begin{lemma}[Two-color tableau transfer]\label{lem:transfer-commute}
For every skew diagram $D$ and all indeterminates $x,y$,
\begin{equation}\label{eq:commute}
  T_xT_y=T_yT_x.
\end{equation}
\end{lemma}
\begin{proof}
The $(I,J)$ entry of $T_xT_y$ is
\[
  \sum_{I\subseteq K\subseteq J}
      x^{\operatorname{col}(K\setminus I)}y^{\operatorname{col}(J\setminus K)}.
\]
Color the cells of the skew shape $J\setminus I$ by $1$ on
$K\setminus I$ and by $2$ on $J\setminus K$.  Since $K$ is an
ideal between $I$ and $J$, this is exactly a reverse plane
partition with entries in $\{1,2\}$ (with the corresponding
coordinate order), and every such partition determines a unique
$K$.  The exponents above are the numbers of columns containing
$1$ and $2$, respectively, i.e., the two components of the
\emph{irredundant content} $\mathrm{ircont}$ used in the definition
of the skew dual stable Grothendieck generating function
$g_{\lambda/\mu}$.  The two-color result needed here is directly the published
Lemma~7 of Galashin--Grinberg--Liu
\cite[Lemma~7]{GalashinGrinbergLiu2016}: for reverse plane partitions
with entries in $\{1,2\}$, it supplies an involution preserving the
column-equality statistic and exchanging the number of columns
containing $1$ and the number containing $2$. Applied to the
intermediate-ideal colorings above, this involution interchanges the
powers of $x$ and $y$ term by term, proving entrywise commutation.
For context, their Theorem~5 states the more general symmetry of the
refined dual stable Grothendieck function, Proposition~4(a) gives its
specialization to $g_{\lambda/\mu}$, and Theorem~6 gives generalized
Bender--Knuth involutions; Section~7 provides a second approach to
the two-color case. The companion Lean development also proves the
two-color identity internally.
\end{proof}

\begin{definition}[Vertical-strip coefficient]\label{def:strip}
A vertical strip has at most one cell in each row.  For a strip $S$,
let $c(S)$ be its number of connected runs.  Define
\[
 q_x(S)=\begin{cases}
 (-x)^{c(S)}(1-x)^{|S|-c(S)},&S\text{ a vertical strip},\\
 0,&\text{otherwise},
 \end{cases}\qquad q_x(\varnothing)=1.
\]
Let $Q_x[I,J]=q_x(J\setminus I)$ for $I\subseteq J$, and zero
otherwise.
\end{definition}

\begin{majorresult}
\begin{lemma}[Inverse and logarithmic derivative]\label{lem:inverse}
For every finite skew diagram,
\begin{equation}\label{eq:inverse}
  Q_xT_x=\mathrm{Id},\qquad
  Q_xT_x'=A_{1-x}.
\end{equation}
In particular $Q_x=T_x^{-1}$, so
\begin{equation}\label{eq:logder}
 A_{1-x}=T_x^{-1}T_x'.
\end{equation}
\end{lemma}
\end{majorresult}
\begin{proof}
Fix $I,J\in\JJ(P_D)$.  If $I\not\subseteq J$, both entries
vanish; if $I=J$, the entries of $Q_xT_x$ and $Q_xT_x'$ are
$1$ and $0$, respectively.  Assume $I\subsetneq J$ and put
$B=J\setminus I$.

\smallskip\noindent\emph{Admissible strips.}
In the matrix product, an intermediate ideal $K$ determines
$S=K\setminus I\subseteq B$.  Because $K$ is a lower ideal,
$S$ is a lower ideal of the induced order on $B$.  A term
$q_x(S)$ can be nonzero only when $S$ is a vertical strip,
i.e., it contains at most one cell in each row.  Therefore every
selected cell must be the leftmost cell of its row in $B$.
Write $[\ell_i,r_i]$ for the occupied row intervals of $B$,
numbered from top to bottom, so that the left endpoints
$\ell_i$ weakly decrease.  If $(i,\ell_i)$ is selected,
any predecessor $(i',j)\in B$ in an earlier row has
$i'<i$ and $j\le\ell_i$; but then
$\ell_{i'}\le j\le\ell_i\le\ell_{i'}$.
Thus $j=\ell_i=\ell_{i'}$: only leftmost cells in earlier
rows with the same left endpoint are forced into $S$.
It follows that, within each maximal consecutive group of rows
having a common left endpoint, the chosen rows form an initial
segment \emph{from the top downward}.  Conversely these choices
are precisely the admissible strips.  The groups are independent.

If a group has $r$ rows, a nonempty chosen initial segment of
length $k$ consists of one connected vertical run and has
weight $-x(1-x)^{k-1}$; the empty choice has weight $1$.
Different groups use different columns and do not meet.
Consequently, if the group lengths are $r_1,\ldots,r_s$ and
$h=\sum_j r_j$, the total strip weight is
\begin{equation}\label{eq:strip-total}
 \sum_{S}q_x(S)
 =\prod_{j=1}^s\left(1-\sum_{k=1}^{r_j}
                x(1-x)^{k-1}\right)
 =(1-x)^h.
\end{equation}

\smallskip\noindent\emph{Removal of the first column.}
Suppose now that $B$ is connected, with
$v=\operatorname{col}(B)$ columns.  Its leftmost column is
$\ell_h$, the first column of the bottom row; it consists exactly
of the leftmost cells in the final group of equal left endpoints.
Every other column retains a cell not chosen by an admissible
strip: this is immediate for columns that are not row-left
endpoints.  If a column $c$ is the left endpoint of an earlier
group, consecutive-row overlap at the transition to the next
group gives a cell in column $c$ that is \emph{not} leftmost
in its lower row.  (Here the overlap inequality is
$\ell_{i+1}<c\le r_{i+1}$.)
Thus $\operatorname{col}(B\setminus S)=v-1$ exactly when $S$
contains the entire last group; otherwise it equals $v$.
The total weight of these exceptional strips is
\begin{equation}\label{eq:strip-special}
 \sum_{S:\,\operatorname{col}(B\setminus S)=v-1}q_x(S)
 =\left(\prod_{j<s}(1-x)^{r_j}\right)
    [-x(1-x)^{r_s-1}]
 =-x(1-x)^{h-1}.
\end{equation}

\smallskip\noindent\emph{Matrix entries.}
For connected nonempty $B$, split the sum according to these
two possible column counts.  Equations
\eqref{eq:strip-total}--\eqref{eq:strip-special} give
\begin{align*}
 [Q_xT_x]_{I,J}
 &=x^v\bigl[(1-x)^h+x(1-x)^{h-1}\bigr]
     +x^{v-1}[-x(1-x)^{h-1}]=0,\\
 [Q_xT_x']_{I,J}
 &=v x^{v-1}\bigl[(1-x)^h+x(1-x)^{h-1}\bigr]
     +(v-1)x^{v-2}[-x(1-x)^{h-1}]\\
 &=x^{v-1}(1-x)^{h-1}=A_{1-x}[I,J].
\end{align*}
When $v=1$, the term involving $(v-1)x^{v-2}$ is simply
the derivative of $x^0=1$, hence zero; no negative power is
needed.  For disconnected $B$, its components occupy disjoint
row and column sets, so the strip weights and $T_x$ entries
factor over components.  The inverse entry contains a vanishing
connected off-diagonal factor, while each differentiated term
contains such a factor when there are at least two components.
Both entries therefore vanish, as required.
\end{proof}

\begin{proposition}[Commuting connected-block kernels]\label{prop:Acommute}
For all $z,u$,
\begin{equation}\label{eq:Acommute}
 A_zA_u=A_uA_z.
\end{equation}
\end{proposition}
\begin{proof}
The polynomial identity \eqref{eq:commute} holds in independent
variables.  Differentiating with respect to each variable yields
commutation among $T_x,T_x',T_y,T_y'$.  Since the inverses belong
to the same commuting algebra, the logarithmic derivatives
$T_x^{-1}T_x'$ also commute.  Use \cref{lem:inverse}.
\end{proof}

\section{Skew-shape block identities}\label{sec:main}
Let $n=|D|$.  The skew-shape identity will be stated for the
corrected statistic of \cref{def:blocks}.  Every integral below has
a polynomial integrand with rational coefficients; its value may be
defined using coefficientwise formal primitives over $\QQ$.

\begin{lemma}[Ordered-root integral]\label{lem:ordered}
For a finite skew diagram $D$, under \cref{def:blocks},
\begin{equation}\label{eq:ordered}
 \frac1{n!}\sum_{\sigma:D\xrightarrow\sim[n]}
     t^{\bl_D(\sigma)}
  =\sum_{k=0}^{n} t^k
     \int_{0<z_1<\cdots<z_k<1}
         [A_{z_1}\cdots A_{z_k}]_{\varnothing,D}
           \,dz_1\cdots dz_k.
\end{equation}
The $k=0$ term is $1$ when $D=\varnothing$ and $0$ otherwise.
\end{lemma}
\begin{proof}
Assign independent uniform real labels in $(0,1)$ to the cells
of $D$.  Almost surely all labels are distinct; their relative
ranks give a uniformly random bijective labeling $\sigma$.
Let $\rho(c)$ be the value of the root of the canonical block
containing $c$, as in \cref{lem:median}.

\smallskip\noindent\emph{Orientation of the ideals.}
The local root inequalities established in \cref{lem:median}
read
\[
 \rho(\text{cell below }c)\le\rho(c)
 \le\rho(\text{cell left of }c).
\]
In particular the root values weakly \emph{decrease} when one
moves south or east through $D$.  If $d$ is northwest of $c$
in the coordinatewise cell order, then
$\rho(d)\ge\rho(c)$.  Therefore, for any threshold $z$,
$\{c\in D:\rho(c)>z\}$ is a lower ideal of $P_D$;
$\{c\in D:\rho(c)\le z\}$ is generally an upper filter,
not a lower ideal.

\smallskip\noindent\emph{The canonical chain.}
List the distinct block-root values in \emph{decreasing} order,
$1>w_1>\cdots>w_k>0$, and define
\[
 I_0=\varnothing,\qquad
 I_j=\{c\in D:\rho(c)\ge w_j\}\quad(1\le j\le k).
\]
Then
$\varnothing=I_0\subsetneq I_1\subsetneq\cdots
\subsetneq I_k=D$ is a chain of lower ideals, and
$B_j=I_j\setminus I_{j-1}$ is the connected canonical block
whose root has value $w_j$.
Conversely, take such a chain with connected nonempty
increments, strictly decreasing root values $w_j$, and all
nonroot boundary inequalities from \cref{def:blocks}.
If two distinct increments are related by the skew-poset
order, their positions in this ideal chain force the
corresponding strict root order.  The uniqueness part of
\cref{lem:median}, including its equivalence with the global
block-order convention, identifies this partition with the
canonical one.  Thus no labeling is counted twice.

\smallskip\noindent\emph{The weights.}
Condition on a fixed chain and fixed root values.  By
\cref{lem:boundary}, $B_j$ has $h(B_j)-1$ nonroot
left-row boundary cells, whose labels lie below $w_j$, and
$\operatorname{col}(B_j)-1$ nonroot bottom-column cells, whose
labels lie above $w_j$.  These two sets are disjoint; the
other nonroot labels are unconstrained.  Independence gives
\[
 w_j^{h(B_j)-1}(1-w_j)^{\operatorname{col}(B_j)-1}
  =A_{w_j}[I_{j-1},I_j].
\]
Integrating over $1>w_1>\cdots>w_k>0$ and summing over
all the chains therefore yields
\[
 \int_{1>w_1>\cdots>w_k>0}
 [A_{w_1}\cdots A_{w_k}]_{\varnothing,D}
     \,dw_1\cdots dw_k.
\]
After putting $z_j=w_{k+1-j}$, the product has factors in
the order $A_{z_k}\cdots A_{z_1}$.  By
\cref{prop:Acommute}, these factors commute and may be
reordered as $A_{z_1}\cdots A_{z_k}$, which is the integrand
in \eqref{eq:ordered}.  Summing over $k$ gives the expected
value of $t^{\bl_D(\sigma)}$, equal to the left side.
For $D=\varnothing$ the empty chain contributes $1$;
otherwise the $k=0$ contribution is zero.
The Lean formalization uses finite label-gap counts and
coefficientwise rational polynomial primitives instead of
continuous probability; see \cref{sec:formalization}.
\end{proof}

\par\smallskip\noindent\emph{Two-cell orientation check.}\\[-3pt]
For two vertically adjacent cells, with upper label $x$ and
lower label $y$, one block occurs when $x<y$, while $x>y$
gives two singleton blocks.  In the latter case the chain is
$\varnothing\subset\{\text{upper cell}\}\subset D$:
the larger root comes first.  Formula \eqref{eq:ordered} gives
$\int_0^1 z\,dz=\tfrac12$ for the one-block contribution
and $\tfrac12$ for the two-block contribution, hence
$\Omega(P_D;t)=(t+t^2)/2$, as required.
\par\smallskip

\begin{majorresult}
\begin{theorem}[Skew-shape block-statistic identity --- Kahane's Conjecture~4.7 (v2)]\label{thm:kahane}
For every finite skew Young diagram $D$, including disconnected
shapes, and with the statistic of \cref{def:blocks},
\begin{equation}\label{eq:kahane}
 n!\,\Om(P_D;t)=
  \sum_{\sigma:D\xrightarrow\sim[n]}t^{\bl_D(\sigma)}.
\end{equation}
This establishes Kahane's Conjecture~4.7 \cite{Kahane2026}
(version~2), with the explicit root-exclusion convention
of \cref{def:blocks}.
\end{theorem}
\end{majorresult}
\begin{proof}
We give a finite polynomial-matrix argument, avoiding a formal
matrix logarithm.  By \cref{lem:inverse}, put
$V(z)=ZQ_{1-z}=ZT_{1-z}^{-1}$.  Its boundary values are
$V(0)=\mathrm{Id}$ and $V(1)=Z$.  Differentiating the inverse
identity and using the commutations in \cref{lem:transfer-commute,prop:Acommute} gives
\begin{equation}\label{eq:V-ode}
 V'(z)=V(z)A_z,\qquad V(0)=\mathrm{Id}.
\end{equation}
Consequently, for each integer $m\ge0$,
$(V(z)^m)'=mV(z)^mA_z$ and $V(0)^m=\mathrm{Id}$.

Define the polynomial-matrix ordered integrals
\[
 R_0(z)=\mathrm{Id},\qquad
 R_k(z)=\int_{0<u_1<\cdots<u_k<z}
      A_{u_1}\cdots A_{u_k}\,du_1\cdots du_k.
\]
They satisfy $R_{k+1}'(z)=R_k(z)A_z$ and
$R_k(0)=0$ for $k>0$.  Because every nonzero entry of $A_z$
moves between \emph{strictly} nested ideals, a product of more than
$n$ such kernels vanishes.  Hence $R_k=0$ for $k>n$, and
\[
 S_m(z)=\sum_{k=0}^{n}m^kR_k(z)
\]
satisfies the same equation and initial condition as $V(z)^m$:
$S_m'=mS_mA_z$, $S_m(0)=\mathrm{Id}$.
These solutions are equal: for their difference, induct on the
cardinality of the target ideal in each matrix entry.  The
strict-upper-triangularity of $A_z$ makes every entry derivative
depend only on previously vanishing entries, and the initial
condition forces its constant term to vanish.  Thus
$S_m(1)=V(1)^m=Z^m$.

Taking the $(\varnothing,D)$ entry of the preceding finite sum, using
\cref{lem:ordered}, and then \cref{lem:zeta}, we get
\[
 \left(\frac1{n!}\sum_\sigma t^{\bl_D(\sigma)}\right)_{t=m}
   =(Z^m)[\varnothing,D]=\Om(P_D;m)
\]
for every positive integer $m$ (and also for $m=0$ by the
empty-shape convention).  Both sides are polynomials in $t$ with
rational coefficients.  Equality at infinitely many arguments proves
\eqref{eq:kahane}.
\end{proof}

\begin{majorresult}
\begin{corollary}[Bivariate consequence of the block identity]\label{thm:bivariate}
As a direct consequence of \cref{thm:kahane}, let $D$ be a skew
diagram with $n$ cells and let $w:P_D\to\ZZ$.  Then
\begin{equation}\label{eq:bivariate}
 n!\,\HH_{P_D,w}(t;q)=
  \sum_{I\in\JJ(P_D)}\frac{n!}{(n-|I|)!}q^{w(I)}
     \sum_{\sigma:P_D\setminus I\xrightarrow\sim[n-|I|]}
       (t-1)^{\bl_{P_D\setminus I}(\sigma)}.
\end{equation}
The empty complement contributes one empty labeling with zero
blocks.  In the variable $u=t-1$, multiplication by a suitable
power of $q$ produces a bivariate polynomial with nonnegative
integer coefficients.
\end{corollary}
\end{majorresult}
\begin{proof}
For each ideal $I$, the complement $P_D\setminus I$ is an induced
filter and the cell poset of a skew diagram, possibly disconnected.
Apply \cref{thm:kahane} to that diagram with $u=t-1$:
\[
 (n-|I|)!\,\Om(P_D\setminus I;u)
   =\sum_{\sigma}u^{\bl_{P_D\setminus I}(\sigma)}.
\]
Substitute in \eqref{eq:intro-H} and multiply by $n!$.  The
factor $n!/(n-|I|)!$ is a nonnegative integer.  Each term on the
right of \eqref{eq:bivariate} is a monomial $q^{w(I)}u^k$
with nonnegative integral multiplicity.  Since the support in the
$q$-direction is finite, a common integral exponent shift makes all
its exponents nonnegative.
\end{proof}

\begin{majorresult}
\begin{corollary}[Two specializations and the common surface subfamily]
\label{cor:common}
For every skew diagram $D$,
\begin{equation}\label{eq:cor-slices}
 \HH_{P_D,w}(t;1)=\frac1{n!}
    \sum_{\sigma:D\xrightarrow\sim[n]}t^{\bl_D(\sigma)},
 \qquad
 \HH_{P_D,w}(2;q)=\sum_{I\in\JJ(P_D)}q^{w(I)}.
\end{equation}
In particular,
\[
  |\JJ(P_D)|=\frac1{n!}
              \sum_{\sigma:D\xrightarrow\sim[n]}2^{\bl_D(\sigma)}.
\]
Whenever the principal-coefficient ideal poset of a plain surface
arc is order-isomorphic to a skew ribbon, its surface-weighted
specialization and its block-statistic specialization are both
computed on that same ideal lattice, subject to the geometric
identification in \eqref{eq:surface-slice}.
\end{corollary}
\end{majorresult}
\begin{proof}
Combine \cref{prop:slices,thm:kahane} and evaluate at $t=2$
for the stated ideal-count identity.  The surface assertion follows from
\cref{lem:fence} and \eqref{eq:surface-slice}, after transporting
the weighting across the order isomorphism.
\end{proof}

\section{Scope and further interpretation}\label{sec:scope}
The bivariate formula \eqref{eq:bivariate} is a consequence of the
skew-shape block identity, not a logically independent proof of a
fourth conjecture.  It offers a common positive expansion for
weighted ideal data and order-polynomial data, without asserting
unimodality of arbitrary weighted specializations.

\medskip\noindent\textbf{Comparison with Kahane's v2.}\quad
We use Kahane's \emph{arXiv:2607.11225v2} (2 November 2026):
Definition~3.2 orders blocks by their root labels;
Proposition~4.4 identifies the leftmost cell of a connected skew
shape; Definition~4.5 specifies the blocks; Lemma~4.6 addresses
their existence and uniqueness; and Conjecture~4.7 is precisely the
skew-shape identity treated in \cref{thm:kahane}.  The block
definition appeared as Definition~4.2 in v1, which explains that
number in Kahane's correspondence of 10 November 2026.  The
correspondence confirms that the comparisons at the root should
be weak, equivalently that the root is omitted from strict
boundary conditions.

For clarity about orientation, Kahane's left parent of $(i,j)$ is
$(i,j-1)$ and his left child is $(i+1,j)$.  Thus the root is the
southwestern boundary cell in our row-coordinate convention.  Our
\cref{lem:boundary,lem:median} expand the associated geometric
properties: the relevant southwest unit square gives
$a\le\rho(s)\le b$, canonical fibers are connected skew subdiagrams,
and local inequalities imply the global block order.
The comparison is also encoded in Lean by
\nolinkurl{CanonicalBlocks.rootLocallyOrdered_iff_kahaneGlobal},
\nolinkurl{admissible_iff_kahaneOrdered}, and
\nolinkurl{correctedKahaneBlocks_exists_unique}.  Our median
construction treats every placement relative to the two neighboring
roots, including the lower-attachment case.  It supplies a direct constructive proof for the
root-excluded convention.

\medskip\noindent\textbf{Geometric and formal scope.}\quad
The annular claims combine an algebraic mutation calculation with a
separate geometric identification of a single bridging curve's
shear-coordinate row in the same triangulation.  The former is
checked in Lean; the latter uses the surface description of
Felikson--Tumarkin and the conventions of Kang--Lee--Lim, and is
not itself represented by a Lean theorem.  Discussions with
A.~Felikson and W.~Kang clarified the shear-sign and coefficient
specialization conventions, while correspondence with Y.~Kahane
clarified the root convention.  These clarifications are not
presented as independent verification or endorsement of the
results.  The intended meaning of a cited geometric object and
its faithful representation in a formal definition remain
mathematical obligations separate from kernel checking.

The general block identity is formalized for diagrams satisfying
\texttt{SkewRowIntervals}, including disconnected shapes, and is
proved at every natural evaluation.  The ordered-root proof uses
finite polynomial primitives and exact label-gap counting rather
than measure theory.  The explicit statement-to-Lean map appears
in \cref{sec:formalization}.  Earlier work on $2\times n$
rectangles by White \cite{White2026} is acknowledged as a special
case and does not enter the general proof.

\section{Lean 4 formalization and statement correspondence}
\label{sec:formalization}
A companion Lean~4 development uses Mathlib and is maintained in
\begin{center}
 \url{https://github.com/adamtilepe/Annular_Obstructions_and_Skew-Shape_Block_Identities}.
\end{center}
The primary source is
\path{Annular_Counterexamples_and_Skew-Shape_Block_Identities.lean}.
The file encodes the order ideal constructions, the algebraic
annular mutations, the corrected canonical block statistic,
tableau-transfer matrices and a finite polynomial approach to the
ordered-root integrals.  In particular, the formalization does not
rely on a numerical test for the general skew-shape identity.

\subsection*{Statements and exact declarations}
The following names are relative to the namespace
\texttt{OrderIdealSpecializations}.  They are identifiers in the
public Lean source, rather than merely informal correspondences.
\begin{center}
\small
\renewcommand{\arraystretch}{1.17}
\begin{tabular}{@{}p{0.33\linewidth}p{0.62\linewidth}@{}}
\toprule
Manuscript statement & Lean declaration(s)\\
\midrule
\cref{prop:slices}, equation~\eqref{eq:slices}
 & \nolinkurl{FinitePoset.H_two},
   \nolinkurl{FinitePoset.H_one_eq_omega_succ}\\
\cref{lem:zeta}, equation~\eqref{eq:zeta}
 & \nolinkurl{FinitePoset.idealZeta_pow_bottom_top_eq_omega}\\
\cref{thm:annulus}, polynomial obstruction
 & \nolinkurl{AnnularFirstMutation.firstPolynomial_not_unimodal},
   \nolinkurl{AnnularFirstMutation.firstPolynomial_two_not_logConcave}\\
\cref{prop:positive}, positive-coefficient obstruction
 & \nolinkurl{AnnularFourthMutation.D_not_logConcave}\\
\cref{lem:transfer-commute}, equation~\eqref{eq:commute}
 & \nolinkurl{TableauTransfer.tableau_transfer_matrices_commute}\\
\cref{lem:median}, bounds and fiber geometry
 & \nolinkurl{CanonicalBlocks.predecessor_bounds_ordered},
   \nolinkurl{CanonicalBlocks.canonicalRootFiber_connectedSkewRows}\\
\cref{lem:inverse}, equation~\eqref{eq:inverse}
 & \nolinkurl{TableauTransfer.vertical_strip_inverse_and_log_derivative}\\
\cref{prop:Acommute}, equation~\eqref{eq:Acommute}
 & \nolinkurl{TableauTransfer.connected_block_kernels_commute}\\
\cref{lem:ordered}, equation~\eqref{eq:ordered}
 & \nolinkurl{OrderedRootIntegral.ordered_root_integral_identity}\\
\cref{thm:kahane}, equation~\eqref{eq:kahane}
 & \nolinkurl{OrderedRootIntegral.skew_shape_block_statistic_identity}\\
\cref{thm:bivariate}, equation~\eqref{eq:bivariate}
 & \nolinkurl{OrderedRootIntegral.weighted_block_refinement_identity}\\
\cref{cor:common}, equation~\eqref{eq:cor-slices}
 & \nolinkurl{OrderedRootIntegral.weighted_block_two_specializations}\\
\bottomrule
\end{tabular}
\end{center}

\subsection*{What is proved and what is not inferred}
The Lean result
\texttt{OrderedRootIntegral.skew\_shape\_block\_statistic\_identity}
assumes \texttt{SkewRowIntervals D}, an explicit representation
of a skew Young diagram by weakly decreasing row endpoints.
Its conclusion is equality at each natural-number evaluation between
the normalized block-generating \emph{polynomial} and the actual
count of isotone colorings.  Since both sides define polynomials,
the familiar polynomial form of \eqref{eq:kahane} follows.

The formalization of \cref{lem:ordered} uses coefficientwise
rational polynomial primitives in place of analytic measure
integrals; the formalization of \cref{thm:kahane} uses
strict-upper-triangular polynomial matrices and a finite differential
equation rather than an abstract exponential and logarithm.
The bivariate result is stated through finite natural-number
evaluations with signed integer $q$-weights; a finite monomial
expansion proves the claimed coefficient nonnegativity.

The annular Lean sections certify their exchange polynomials and
coefficient inequalities; they do \emph{not} construct a topological
annulus or prove that a frozen row is the shear vector of a curve.
The block-statistic theorem \cref{thm:kahane} certifies its precisely
stated root-excluded identity, which is compared with Kahane's
Conjecture~4.7 (arXiv:2607.11225v2) in \cref{sec:scope}.
External geometric identifications are not assertions of the Lean
kernel.

\subsection*{Reproducibility and axiom audit}
The repository records its Lean toolchain and Mathlib revision.
From the project root, the compilation and principal audit may be run
using
\begin{center}
\texttt{lake env lean Annular\_Counterexamples\_and\_Skew-Shape\_Block\_Identities.lean}
\end{center}
The source includes the command
\begin{center}
\small\texttt{\#print axioms}\quad
\nolinkurl{OrderIdealSpecializations.OrderedRootIntegral.skew_shape_block_statistic_identity}
\end{center}
The successful Lean run reported precisely the standard axioms
\texttt{propext}, \texttt{Classical.choice}, and \texttt{Quot.sound},
with no additional project-specific axiom.  These are standard
foundational axioms of Lean/Mathlib.  As always, formal kernel
checking establishes the encoded statement, while agreement of the
encoding with an external mathematical description is a separate
question.


\begin{thebibliography}{GGL16}
\bibitem[Kah26]{Kahane2026}
Y.~Kahane, \emph{Combinatorial interpretation of the coefficients of the
order polynomial of fence posets}, arXiv:2607.11225v2 (2026).
\url{https://arxiv.org/abs/2607.11225}.

\bibitem[KLL26]{KangLeeLim2026}
W.~Kang, K.~Lee and E.~Lim,
\emph{Unimodality and cluster algebras from surfaces},
European J. Combin. \textbf{136} (2026), 104392.
\href{https://doi.org/10.1016/j.ejc.2026.104392}{doi:10.1016/j.ejc.2026.104392}.

\bibitem[FT24]{FeliksonTumarkin2024}
A.~Felikson and P.~Tumarkin,
\emph{Cluster algebras of finite mutation type with coefficients},
J. Comb. Algebra \textbf{8} (2024), 375--418.
\href{https://doi.org/10.4171/JCA/92}{doi:10.4171/JCA/92}.

\bibitem[GGL16]{GalashinGrinbergLiu2016}
P.~Galashin, D.~Grinberg and G.~Liu,
\emph{Refined dual stable Grothendieck polynomials and generalized
Bender--Knuth involutions}, Electron. J. Combin.
\textbf{23}(3) (2016), P3.14.
\href{https://doi.org/10.37236/5737}{doi:10.37236/5737}.

\bibitem[LP07]{LamPylyavskyy2007}
T.~Lam and P.~Pylyavskyy,
\emph{Combinatorial Hopf algebras and $K$-homology of Grassmannians},
Int. Math. Res. Not. \textbf{2007} (2007), rnm125.
\href{https://doi.org/10.1093/imrn/rnm125}{doi:10.1093/imrn/rnm125}.

\bibitem[FZ07]{FominZelevinsky2007}
S.~Fomin and A.~Zelevinsky,
\emph{Cluster algebras IV: Coefficients},
Compos. Math. \textbf{143}(1) (2007), 112--164.
\href{https://doi.org/10.1112/S0010437X06002521}{doi:10.1112/S0010437X06002521}.

\bibitem[MSW11]{MusikerSchifflerWilliams2011}
G.~Musiker, R.~Schiffler and L.~Williams,
\emph{Positivity for cluster algebras from surfaces},
Adv. Math. \textbf{227}(6) (2011), 2241--2308.
\href{https://doi.org/10.1016/j.aim.2011.04.018}{doi:10.1016/j.aim.2011.04.018}.

\bibitem[Whi26]{White2026}
P.~White, \emph{A closed form for the block-statistic generating function
of the $2\times n$ rectangle: the two-row case of Kahane's conjecture},
arXiv:2609.17561 (2026).
\url{https://arxiv.org/abs/2609.17561}.
\end{thebibliography}
\end{document}